\documentclass[a4paper,12pt]{amsart}

\usepackage{amsmath}
\usepackage{amssymb}
\usepackage{mathrsfs}
\usepackage{ifthen}
\usepackage{graphicx}
\usepackage[T1]{fontenc} 

\def\switchlinenumbers{\@ifstar
    {\let\makeLineNumberOdd\makeLineNumberRight
     \let\makeLineNumberEven\makeLineNumberLeft}%
    {\let\makeLineNumberOdd\makeLineNumberLeft
     \let\makeLineNumberEven\makeLineNumberRight}%
    }

\def\setmakelinenumbers#1{\@ifstar
  {\let\makeLineNumberRunning#1%
   \let\makeLineNumberOdd#1%
   \let\makeLineNumberEven#1}%
  {\ifx\c@linenumber\c@runninglinenumber
      \let\makeLineNumberRunning#1%
   \else
      \let\makeLineNumberOdd#1%
      \let\makeLineNumberEven#1%
   \fi}%
  }

\nonstopmode \numberwithin{equation}{section}
\newtheorem*{theorem*}{Theorem}

\newtheorem{thm}{Theorem}[section]
\newtheorem{cor}[equation]{Corollary}
\newtheorem{lem}[equation]{Lemma}
\newtheorem{prop}[equation]{Proposition}

\theoremstyle{definition}
\newtheorem{defn}{Definition}[section]

\newtheorem{prob}[equation]{Problem}
\newtheorem{rem}{Remark}[section]

\newtheorem{innercustomthm}{Theorem}
\newenvironment{customthm}[1]
  {\renewcommand\theinnercustomthm{#1}\innercustomthm}
  {\endinnercustomthm}

\newcounter{minutes}
\divide\time by 60
\newcounter{hours}
\multiply\time by 60
\addtocounter{minutes}{-\time}

\newcounter {own}
\def\theown {\thesection       .\arabic{own}}

\newenvironment{pf}[1][]{%
 \vskip 3mm
 \noindent
 \ifthenelse{\equal{#1}{}}%
  {{\slshape Proof. }}%
  {{\slshape #1.} }%
 }%
{\qed\bigskip}

\newcounter{alphabet}

\def\be{\begin{equation}}
\def\ee{\end{equation}}

\newcommand{\bee}{\begin{enumerate}}
\newcommand{\eee}{\end{enumerate}}

\newcommand{\blem}{\begin{lem}}
\newcommand{\elem}{\end{lem}}
\newcommand{\bthm}{\begin{thm}}
\newcommand{\ethm}{\end{thm}}
\newcommand{\bcor}{\begin{cor}}
\newcommand{\ecor}{\end{cor}}
\newcommand{\beg}{\begin{examp}}
\newcommand{\eeg}{\end{examp}}
\newcommand{\begs}{\begin{examples}}
\newcommand{\eegs}{\end{examples}}
\newcommand{\bdefe}{\begin{defin}}
\newcommand{\edefe}{\end{defin}}
\newcommand{\bprob}{\begin{prob}}
\newcommand{\eprob}{\end{prob}}
\newcommand{\bei}{\begin{itemize}}
\newcommand{\eei}{\end{itemize}}

\newcommand{\norm}[1]{\left\lVert#1\right\rVert}

\begin{document}

\title{On multidimensional Bohr radii for Banach spaces}

\author{Vasudevarao Allu}
\address{Vasudevarao Allu,
Department of Mathematics,
School of Basic Sciences,
Indian Institute of Technology Bhubaneswar,
Bhubaneswar-752050, Odisha, India.}
\email{avrao@iitbbs.ac.in}


\author{Subhadip Pal}
\address{Subhadip Pal,
	Department of Mathematics,
	School of Basic Sciences,
	Indian Institute of Technology Bhubaneswar,
	Bhubaneswar-752050, Odisha, India.}
\email{subhadippal33@gmail.com}

\subjclass[{AMS} Subject Classification:]{Primary 32A05, 32A10, 46B07; Secondary 46B09, 46G20}

\def\thefootnote{}
\footnotetext{ {\tiny File:~\jobname.tex,
printed: \number\year-\number\month-\number\day,
          \thehours.\ifnum\theminutes<10{0}\fi\theminutes }
} \makeatletter\def\thefootnote{\@arabic\c@footnote}\makeatother

\begin{abstract}
In this paper, we study a more general version of multidimensional Bohr radii for the holomorphic functions defined on unit ball of $\ell^n_q\,\,(1\leq q\leq \infty)$ spaces with values in arbitrary complex Banach spaces. More precisely, we study the multidimensional Bohr radii for bounded linear operators between complex Banach spaces, primarily motivated by the work of A. Defant, M. Maestre, and U. Schwarting [Adv. Math. 231 (2012), pp. 2837--2857]. We obtain the exact asymptotic estimates of multidimensional Bohr radius for both finite and infinite dimensional Banach spaces. As an application, we find the lower bound of arithmetic Bohr radius.
\end{abstract}

\maketitle
\pagestyle{myheadings}
\markboth{Vasudevarao Allu and Subhadip Pal}{On multidimensional Bohr radii for Banach spaces}

\section{Introduction}
 In the beginning of the 20th century while studying the well-known Riemann $\zeta$-function, a Danish mathematician Harald Bohr made significant developments to the general theory of Dirichlet series. A Dirichlet series is an summation expression of the form
\begin{equation}\label{eqn-001}
	D(s):=\sum_{n\geq 1}\frac{a_n}{n^s},
\end{equation}
where $a_n\in \mathbb{C}$ and $s=\sigma +it$ is a complex variable. It is well-known that the region of convergence, absolutely convergence, and uniform convergence of a Dirichlet series define to be the half planes of the form $\{s\in \mathbb{C}: \mbox{Re}(s)>\sigma_{0}\}$ in the complex plane $\mathbb{C}$. Bohr was mainly interested in tackling the region of convergence of this Dirichlet series. More precisely, Bohr was investigating the width of the greatest strip for which a Dirichlet series converges uniformly but not absolutely. This problem is known as \textit{absolute convergence problem} of Bohr in the theory of Dirichlet series. In order to tackle this problem, Bohr has obtained a nice connection between a Dirichlet series and a power series of infinite variables. Indeed, given a Dirichlet series of the form \eqref{eqn-001}, Bohr considered for each $n\in \mathbb{N}$ the prime decomposition of $n=p^{\alpha_1}_{1}\cdots p^{\alpha_r}_{r}$ (where $(p_m)_{m\in \mathbb{N}}$ denotes the sequence of order prime numbers) and defined $z=(p^{-s}_{1}, \ldots, p^{-s}_{r})$. Therefore, Bohr has established the following identity:
\begin{equation}\label{eqn-002}
	D(s)=\sum_{n\geq 1}a_n(p^{-s}_{1})^{\alpha_1}\cdots (p^{-s}_{r})^{\alpha_r}=\sum a_n z^{\alpha_1}_{1}\cdots z^{\alpha_r}_r.
\end{equation}
The above correspondence \eqref{eqn-002} is known as the \textit{Bohr transform}. It helps us to translate the problem on Dirichlet series in terms of power series and thus, one can apply complex analysis techniques to solve this problem. As a part of these studies, Bohr once asked whether is it possible to compare the absolute value of a power series in one complex variable with the sum of the absolute values of its coefficients. In 1914, Bohr \cite{Bohr-1914} proved the following remarkable result which nowadays known as \textit{Bohr's inequality}:

\begin{customthm}{A} \cite{Bohr-1914}\label{Bohr-Thm-1914}
	If $f(z)= \sum_{n=0}^{\infty}a_n z^n$ is a holomorphic function in unit disk such that $\sup_{z\in \mathbb{D}}|f(z)|<\infty$, then 
	\begin{align}\label{eqn-003}
		\sum_{n=0}^{\infty}|a_n|r^n \leq \sup_{z\in \mathbb{D}}|f(z)|
	\end{align}
	for $|z|=r \leq 1/3$, and the constant $1/3$, referred to as the Bohr radius, which is sharp.
\end{customthm}

\noindent
In fact, Bohr himself obtained the radius as $1/6$. Later, this radius was further improved to $1/3$ by M. Riesz, I. Schur, and N. Wiener \cite{sidon-1927,tomic-1962}, independently. This significant inequality was forgotten for long eighty years until Dixon \cite{Dixon & BLMS & 1995} used the inequality in his paper to disprove a long-standing conjecture in Operator Algebras: 	If a Banach algebra satisfies the von Neumann inequality then it is necessarily an operator algebra, \textit{i.e.}, which is isometrically isomorphic to a closed subalgebra of $\mathcal{L}(H)$, the algebra of all bounded operators for some Hilbert space $H$. After that, several researchers like Aizenberg \cite{aizn-2000a,aizn-2000b,aizenberg-2001}, Balasubramanian, Calado and Queff\'{e}lec \cite{bala-studia-2006}, Boas and Khavinson \cite{boas-1997}, Boas \cite{boas-2000}, Bombieri and Bourgain \cite{bombieri-2004}, Galicer, Mansilla and Muro \cite{galicer-TAMS-2020}, Hamada, Honda and Kohr \cite{hamada-israel-2009}, Paulsen and Singh \cite{paulsen-2002}, Popescu \cite{popescu-2019} have paid much attention to this result to use it in different contexts of mathematics and generalized it. In this paper, we mostly focus on the unit ball of Minkowski space $\ell^n_q$, where $q \in [1,\infty]$. For $q \in [1,\infty)$, we denote
$$B_{\ell^n _q}= \left\{z \in \mathbb{C}^n:\norm{z}_{q}=\left(\sum_{i=1}^{n}|z_{i}|^q\right)^{1/q}<1\right\},$$ 
and 
$$B_{\ell^n _\infty}= \left\{z \in \mathbb{C}^n: \norm{z}_{\infty}=\sup _{1\leq i \leq n} |z_{i}|<1\right\}.$$

Before proceeding for further discussions, let us first introduce the following notion of multidimensional Bohr radius for holomorphic functions defined on $B_{\ell^n_q},\, 1\leq q\leq \infty$ with values in a complex Banach space $X$, primarily motivated by the work of Defant {\it et al.} \cite{defant-adv-math-2012}.
\begin{defn}\label{Pal-Vasu-P5-defn-01}
	Let $1\leq q\leq \infty$, $1\leq p<\infty$, $\lambda\geq 1$, $n\in \mathbb{N}$, and let $T: X \rightarrow Y$ be a bounded linear operator between complex Banach spaces $X$ and $Y$ such that $\norm{T}\leq \lambda^{1/p}$. For $1\leq p<\infty$, the $\lambda_{p}$-Bohr radius of $T$, denoted by $K^p(B_{\ell^n_{q}}, T,\lambda)$, is the supremum of all $r\geq 0$ such that for all holomorphic functions $f(z)=\sum_{\alpha\in \mathbb{N}^n_{0}}x_{\alpha}z^{\alpha}$ on $B_{\ell^n_{q}}$ we have 
	\begin{equation}\label{Pal-Vasu-P5-eqn-int-04}
		\sup_{z\in rB_{\ell^n_{q}}}\sum_{\alpha\in \mathbb{N}^n_{0}}\norm{T(x_{\alpha})z^{\alpha}}^{p}_{Y}\leq \lambda\sup_{z\in B_{\ell^n_{q}}}\norm{\sum_{\alpha\in \mathbb{N}^n_{0}}x_{\alpha}z^{\alpha}}^p_{X}.
	\end{equation}
\end{defn}
\noindent
We write $K^p(B_{\ell^n_q},T)$ for $\lambda=1$. If $T$ is the identity map on $X$ then we use the notation $K^p(B_{\ell^n_q},X,\lambda)$ and $K^p(B_{\ell^n_q},X)$. For $X=\mathbb{C}$, we write $K^p(B_{\ell^n_q},\lambda)$ for $K^p(B_{\ell^n_q},\mathbb{C},\lambda)$ and $K(B_{\ell^n_q})$ whenever $p=1$ and $\lambda=1$.
\vspace{3mm}

We observe that if $f$ is an unbounded holomorphic function on $B_{\ell^n_q}$ then the inequality \eqref{Pal-Vasu-P5-eqn-int-04} holds trivially. Therefore, throughout our discussion, we assume the holomorphic functions to be bounded on $B_{\ell^n_q}$. Clearly, with the above notation, Theorem \ref{Bohr-Thm-1914} reduces to $	K(\mathbb{D})=1/3$.

In the classical sense, Bombieri and Bourgain \cite{bombieri-1962,bombieri-2004} have studied the constant $K(\mathbb{D},\lambda)$ whenever $\lambda\geq 1$. In fact, Bombieri \cite{bombieri-1962} has obtained the exact value 
\begin{equation*}
	K(\mathbb{D},\lambda)=\frac{1}{3\lambda-2\sqrt{2(\lambda^2-1)}},
\end{equation*}
 where $\lambda \in [1,\sqrt{2}]$, and Bombieri and Bourgain \cite{bombieri-2004} have determined the exact asymptotic behavior of $K(\mathbb{D},\lambda)$ as $\sqrt{\lambda^2-1}/\lambda$ whenever $\lambda \rightarrow \infty$. In 2010, Blasco \cite{Blasco-OTAA-2010} established that $K(\ell^2_q)=0$ for $1\leq q\leq \infty$. Thus, the study of Bohr radius problem becomes irrelevant for $\text{dim}(X)>1$. Hence, it seems reasonable to modify the Bohr inequality to obtain more interesting result. That is why the authors \cite{defant-adv-math-2012} have introduced $\lambda$ in the definition of Bohr radius $K(B_{\ell^n_{\infty}},T,\lambda)$ of bounded linear operator $T$. In \cite{defant-adv-math-2012}, the authors have established that $K(B_{\ell^n_{\infty}},T,\lambda)>0$ for every bounded linear operator $T:X\rightarrow Y$ between complex Banach spaces $X$ and $Y$. Later, in this direction, the work was extended by Kumar {\it et al.} \cite{kumar-2023-arxiv} and studied the Bohr radius $K(B_{\ell^n_q},T,\lambda)$ for any bounded operator $T$, $\lambda\geq 1$, and $1\leq q\leq \infty$.\\[0.5mm]
 
 In $1989$, Dineen and Timoney \cite{Dineen-Timoney-1989} extensively studied the constant $K(B_{\ell^n _\infty})$ and later, their result has been quickly refined in \cite{boas-1997}. In $1997$, Boas and Khavinson \cite{boas-1997} proved the following estimate for each $n \in \mathbb{N}$ with $n \geq 2$,
 \begin{equation*}
 	\frac{1}{3\sqrt{n}} \leq K(B_{\ell^n_{\infty}}) \leq 2 \sqrt{\frac{\log n}{n}}.	
 \end{equation*}
 In view of results of Aizenberg, Boas, Dineen, Khavinson and Timoney (see \cite{aizn-2000a},\cite{boas-1997}, \cite{boas-2000}, \cite{Dineen-Timoney-1989}), estimating the multidimensional Bohr radii, for every $1\leq q\leq \infty$ and all $n$ we have
 \begin{equation*}
 	\frac{1}{C}\left(\frac{1}{n}\right)^{1-\frac{1}{\min\{q,2\}}}\leq K(B_{\ell^n_q})\leq C\left(\frac{\log n}{n}\right)^{1-\frac{1}{\min\{q,2\}}},
 \end{equation*}
 where $C$ is a constant independent of $q$ and $n$. One of the main reasons behind the occurrence of $\log$-term in the upper estimates of all known non-trivial results on multidimensional Bohr radii is the probabilistic methods used by the authors. In \cite[p. 326]{boas-2000} Boas conjectured that {\it ``This logarithmic factor, an artifact of the proof, presumably should not really be present''}. In 2006, Defant and Frerick \cite{defant-2006} disproved the conjecture by obtaining a logarithmic lower bound which is almost correct asymptotic estimate for the Bohr radius $K(B_{\ell^n_q})\,\, 1\leq q\leq \infty$. Using the methods developed in \cite{defant-2003}, Defant and Frerick \cite{defant-2006} have proved that there is a constant $C>0$ such that for every $1\leq q\leq \infty$ and $n$ we have
 \begin{equation}\label{Pal-Vasu-P5-int-eqn-05}
 	\frac{1}{C}\left(\frac{\log n/\log \log n}{n}\right)^{1-\frac{1}{\min\{q,2\}}}\leq K(B_{\ell^n_q}).
 \end{equation}
 Later, \eqref{Pal-Vasu-P5-int-eqn-05} was improved by the same authors in \cite{defant-angew-2011}. The case for $q=\infty$ was settled by Defant and his co-authors \cite{defant-2011} and they showed that $K(B_{\ell^n_{\infty}})$ behaves asymptotically as $\sqrt{\log n/n}$ as a consequence of hypercontractivity of Bohnenblust-Hille inequality. Further, Bayart {\it et al.} \cite{bayart-advance-2014} have extended the techniques and obtained the following:
 \begin{equation*}
 	\lim_{n\rightarrow \infty}\frac{K(B_{\ell^n_{\infty}})}{\sqrt{\frac{\log n}{n}}}=1.
 \end{equation*}
 
 The powered versions of the Bohr radii $K^p(\mathbb{D})$ and $K^p(B_{\ell^n_{\infty}})$ for scalar-valued holomorphic functions were first studied by Djakov and Ramanujan \cite{Djakov & Ramanujan & J. Anal & 2000}. Later, B\'{e}n\'{e}teau {\it et al.} \cite{bene-2004} explored the analogous radii in the context of Hardy spaces. Apart from the scalar valued functions, Blasco \cite{Blasco-Collect-2017} has investigated the radius $K^p(\mathbb{D},X)$ for arbitrary complex Banach spaces. In fact, Blasco \cite[Theorem 1.10]{Blasco-Collect-2017} has proved that the Bohr radius $K^p(\mathbb{D},X)>0$ if, and only if, $X$ is $p$-uniformly $\mathbb{C}$-convex, where $2\leq p<\infty$. Recently, Das \cite{das-2023,Das-2024-forum} has showed the precise asymptotic value of $K^p(B_{\ell^n_{\infty}})$ and investigated $p$-Bohr radius for Hardy space of $X$-valued holomorphic functions defined on complete Reinhardt domain.\\[1mm]
 
   To explore more general aspects of multidimensional Bohr radii, we introduce $\lambda$ in our Definition \ref{Pal-Vasu-P5-defn-01}, in more general sense. Indeed, we shall prove that for every $1\leq p<\infty$ and bounded linear operator $T$ with $\norm{T}<\lambda^{1/p}$, the multidimensional Bohr radius $K^p(B_{\ell^n_q},T,\lambda)$ is nonzero. Further, our aim is to study the correct asymptotic estimates for the Bohr radii $K^p(B_{\ell^n_q},T,\lambda)$ for arbitrary Banach spaces. As a consequence, we obtain lower bound of arithmetic Bohr radius for Banach spaces. The organization of this paper is as follows. After the preliminaries, Section \ref{Pal-Vasu-P5-section-03} accumulates all the main results of this paper. In this section, Proposition \ref{Pal-Vasu-P5-prop-01} shows the for every bounded operator $T$ with $\norm{T}<\lambda^{1/p}$ the Bohr radius $K^p(B_{\ell^n_q},T,\lambda)$ is positive. Further, Theorem \ref{Pal-Vasu-P5-thm-01} and Theorem \ref{Pal-Vasu-P5-thm-02} give the asymptotic behavior of $K^p(B_{\ell^n_q},X,\lambda)$ for finite and infinite dimensional Banach space $X$ respectively, in case of identity operator on $X$. Finally, in Theorem \ref{Pal-Vasu-P5-thm-03}, we obtain lower estimates for $K^p(B_{\ell^n_q},T,\lambda)$. Section \ref{Pal-Vasu-P5-section-04} keeps two important lemmas which play vital role to prove the main results. Section \ref{Pal-Vasu-P5-section-05} contains the proof of main results.
 
\section{Preliminaries}
We employ standard terminologies and notions from classical theory of Banach spaces (see e.g. \cite{lindenstrauss-book-I,lindenstrauss-book-II}). We denote $\mathbb{D}$ to be the open unit disk in complex plane $\mathbb{C}$. Throughout the paper, it has been assumed that every Banach space $X$ is complex. $X^*$ will denote the topological dual of Banach space $X$ and we write $B_{X}$ for the open unit ball in $X$. As usual, we denote $\ell^n_q$ for the Banach space of all $n$-tuples $z=(z_1, \ldots, z_n)\in \mathbb{C}^n$ endowed with the norm 
\begin{equation*}
	\norm{z}_q:= \left(\sum_{i=1}^{n}|z_i|^q\right)^{\frac{1}{q}}\,\,\, \mbox{for} \,\, 1\leq q<\infty
\end{equation*}
and $\norm{z}_{\infty}:= \max_{1\leq i\leq n}|z_i|$ for $q=\infty$. For $1\leq q\leq \infty$, the conjugate exponent $q'$ for $q$ is defined by $1/q+1/q'=1$. For $x\in X$, the absolute value of $x$ is defined by $|x|=x\vee (-x)$. Recall that a Banach lattice is a Banach space $X$ which is a vector lattice with $|x|\leq |y|$ implies $\norm{x}\leq \norm{y}$ for all $x,y\in X$. A Banach lattice is said to be $r$-concave, where $1\leq r<\infty$, if there is a constant $M>0$ such that 
\begin{equation}\label{Pal-Vasu-P5-pre-e-01}
	\left(\sum_{j=1}^{n}\norm{x_j}^r\right)^{\frac{1}{r}}\leq M \norm{\left(\sum_{j=1}^{n}|x_j|^r\right)^{\frac{1}{r}}}
\end{equation} 
for every arbitrarily chosen finite elements $x_1,\ldots, x_n\in X$. The best such constant satisfying \eqref{Pal-Vasu-P5-pre-e-01} is as usual denoted by $M_r(X)$. For $2\leq r\leq \infty$, a Banach space $X$ is said to have cotype $r$, if there is a constant $C>0$ such that for every choice of finitely many vectors $x_1, \ldots, x_n\in X$, we have
\begin{equation}\label{Pal-Vasu-P5-pre-e-02}
	\left(\sum_{j=1}^{n}\norm{x_j}^r\right)^{\frac{1}{r}} \leq C \left(\int_{0}^{1}\norm{\sum_{j=1}^{n}r_j(t)x_j}^2\, dt\right)^{\frac{1}{2}},
\end{equation}
where $r_j$ is the $j$th Rademacher function on $[0,1]$. Further, the best such constant $C$ satisfying \eqref{Pal-Vasu-P5-pre-e-02} is usually denoted by $C_r(X)$. We write 
\begin{equation*}
	\text{Cot}(X):= \inf \left\{2\leq r\leq \infty: X \,\, \mbox{has cotype}\, r\right\}.
\end{equation*}
For a Banach space with cotype $\infty$ {\it i.e.}, $\text{Cot}(X)=\infty$, we denote $\frac{1}{\text{Cot}(X)}=0$. It is to be noted that the notions of concavity and cotype are closely related concepts in the context of Banach lattices. In fact, a $r$-concave Banach lattice with $r\geq 2$ is of cotype $r$. On the other way, each Banach lattice of cotype $2$ is $2$-concave, and a Banach lattice of cotype $r>2$ is $s$-concave for all $s>r$.\\[1mm]

Given $m,n\in \mathbb{N}$, we consider the following sets of indices:
\begin{align*}
	\mathcal{M}(m,n)&=\left\{\mathbf{j}=(j_1,\ldots,j_m): 1\leq j_1,\ldots,j_m\leq n\right\}\\ &
	=\left\{1,2,\ldots,n\right\}^m,
\end{align*} 
\begin{equation*}
	\mathcal{J}(m,n)=\left\{\mathbf{j}\in \mathcal{M}(m,n): 1\leq j_1\leq \cdots\leq j_m\leq n\right\}.
\end{equation*}
Furthermore, we write another closely related multi-index set
\begin{equation*}
	\Lambda(m,n)=\left\{\alpha=(\alpha_1,\ldots,\alpha_n)\in \mathbb{N}^n_{0}: |\alpha|=m\right\},
\end{equation*}
where $|\alpha|=\alpha_1+\cdots+\alpha_n$. We may identify $\Lambda(m,n)$ with $	\mathcal{J}(m,n)$ and vice-versa. For $\alpha=(\alpha_1,\ldots,\alpha_n)\in \Lambda(m,n)$, the associated index $\mathbf{j}_{\alpha}\in \mathcal{J}(m,n)$ is given by $\mathbf{j}_{\alpha}=\left(1,\stackrel{\alpha_1}{\ldots},1, 2,\stackrel{\alpha_{2}}{\ldots}2, \ldots,n,\stackrel{\alpha_n}{\ldots},n\right)$.
Moreover, for $\mathbf{j}\in \mathcal{J}(m,n)$ the associated multi-index $\mathbf{j}_{\alpha}\in \Lambda(m,n)$ is given by $\mathbf{j}_{\alpha_r}=|\{k:j_k=r\}|$. This gives a one-to-one correspondence between $\Lambda(m,n)$ and $\mathcal{J}(m,n)$. To identify $\mathcal{M}(m,n)$, we define the following equivalence relation on $\mathcal{M}(m,n)$: for $\mathbf{i}, \mathbf{j}\in \mathcal{M}(m,n)$, we denote $\mathbf{i} \thicksim \mathbf{j}$ if, and only if, there is a permutation $\sigma \in S_m$ such that $\mathbf{j}=(i_{\sigma(1)},\ldots,i_{\sigma(m)})$. We observe that 
\begin{equation*}
	\mathcal{M}(m,n)=\cup_{\mathbf{j}\in \mathcal{J}(m,n)}[\mathbf{j}].
\end{equation*} 
Moreover, it is worth noting that for each $\alpha\in \Lambda(m,n)$ we have
\begin{equation*}
	\text{card}[\mathbf{j}_{\alpha}]=\frac{m!}{\alpha!},
\end{equation*}
where $\alpha!=\alpha_1!\cdots\alpha_n!$. Further, we consider the following notation from \cite{bayart-J-Anal-2019}: given a subset $\mathcal{J}\subset \mathcal{J}(m,n)$, we denote
\begin{equation*}
	\mathcal{J}^*=\left\{\mathbf{j}\in \mathcal{J}(m-1,n): \,\mbox{there is}\,\,k\geq 1, (\mathbf{j},k)\in \mathcal{J}\right\}.
\end{equation*}
\vspace{1mm}

For $1\leq q\leq \infty$, let $\mathcal{P}(^m\ell^n_q,X)$ be the linear space of all $m$-homogeneous polynomial $Q:\ell^n_q\rightarrow X$ of $n$ complex variables with values in $X$, equipped with uniform (or sup) norm 
\begin{equation*}
	\norm{Q}_{\ell^n_q}:=\sup_{z\in B_{\ell^n_{q}}}\norm{Q(z)}_{X}.
\end{equation*}
Using the above multi-index notations, a polynomial $Q\in \mathcal{P}(^m\ell^n_q,X)$
can be expressed in the form 
\begin{equation*}
	Q(z)=\sum_{\alpha\in\Lambda(m,n)}x_{\alpha}z^{\alpha}=\sum_{\mathbf{j}\in \mathcal{J}(m,n)} c_{\mathbf{j}}z_{\mathbf{j}},
\end{equation*}
where $x_{\alpha}, c_{\mathbf{j}}\in X$ and $x_{\alpha}=c_{\mathbf{j}}$ with $\mathbf{j}=\left(1,\stackrel{\alpha_1}{\ldots},1, 2,\stackrel{\alpha_{2}}{\ldots}2, \ldots,n,\stackrel{\alpha_n}{\ldots},n\right)$. If $f:B_{\ell^n_q}\rightarrow X$ is an $X$-valued holomorphic function then we write $f(z)=\sum_{\alpha\in \mathbb{N}^n_{0}}x_{\alpha}(f)z^{\alpha}$ as the monomial series expansion of $f$, where $x_{\alpha}(f)=(\partial^{\alpha}f(0))/\alpha!$ is the $\alpha$th coefficient of the expansion.\\[1mm]

A bounded linear operator $T:X\rightarrow Y$ between two Banach spaces $X$ and $Y$ is called $(r,s)$-summing, $1\leq r,s<\infty$, if there is a constant $\widetilde{C}>0$ such that for each choice of finitely many vectors $x_1,\ldots,x_n\in X$, we have
\begin{equation}\label{Pal-Vasu-P5-pre-e-03}
	\left(\sum_{j=1}^{n}\norm{T(x_j)}^r\right)^{\frac{1}{r}}\leq \widetilde{C}\sup_{x^*\in B_{X^*}}\left(\sum_{j=1}^{n}|x^*(x_j)|^s\right)^{\frac{1}{s}}.
\end{equation}
Also, the best such constant $\widetilde{C}$ satisfying \eqref{Pal-Vasu-P5-pre-e-03} is as usual denoted by $\pi_{r,s}(T)$. In particular, if $r=s$, we call $T$ is $r$-summing and we write the best constant as $\pi_r(T)$.
\section{Main Results}\label{Pal-Vasu-P5-section-03}
As a crucial property, we ensure that for every bounded operator $T$ with $\norm{T}<\lambda^{1/p}$, the $\lambda_{p}$-Bohr radius of $T$ is non-zero.

 \begin{prop}\label{Pal-Vasu-P5-prop-01}
 	Let $1\leq p<\infty$ and $\lambda >1$. Suppose $T:X\rightarrow Y$ is a bounded operator between complex Banach spaces $X$ and $Y$ such that $\norm{T}<\lambda^{1/p}$. Then there exists a constant $M>0$ such that for all $1\leq q\leq \infty$ and $n\in \mathbb{N}$ we have 
 	\begin{equation*}
 		K^p(B_{\ell^n_q},T,\lambda)\geq M \frac{1}{n^{\frac{1}{p}-\frac{1}{pq}}},
 	\end{equation*}
 	where 
 	$$
 	M=\left\{\begin{array}{ll}
 		\max \left\{\left(\frac{\lambda-\norm{T}^{p}}{2\lambda-\norm{T}^p}\right)^{1/p}, \left(\frac{\lambda-\norm{T}^p}{(\lambda-\norm{T}^p+1)\norm{T}}\right)^{1/p}\right\}, & \mbox{ if  \,$\norm{T}\geq 1$},\\[5mm] 
 		\max \left\{\left(\frac{\lambda-\norm{T}^{p}}{2\lambda-\norm{T}^p}\right)^{1/p}, \left(\frac{\lambda-\norm{T}^p}{\lambda-\norm{T}^p+1}\right)^{1/p}\right\}, & \mbox{ if \,$0<\norm{T}<1$}.
 	\end{array}\right.
 	$$
 \end{prop}
As a consequence, whenever we consider the identity operator on $X$, Proposition \ref{Pal-Vasu-P5-prop-01} demonstrates the lower bound for the classical $p$-Bohr radius of Banach spaces.

\begin{cor}\label{Pal-Vasu-P5-cor-01}
	Let $X$ be a complex Banach space and $\lambda>1$, $1\leq p<\infty$, Then for all $1\leq q\leq \infty$ and $n\in \mathbb{N}$ we have
	\begin{equation*}
		K^p(B_{\ell^n_q},X,\lambda)\geq \left(\frac{\lambda-1}{\lambda}\right)^{\frac{1}{p}}\frac{1}{n^{\frac{1}{p}-\frac{1}{pq}}}.
	\end{equation*}
\end{cor}
\begin{rem}
	It is worth noting that for $\lambda=1$ and $p=1$, we actually obtain the trivial lower bound $K(B_{\ell^n_q},X)\geq 0$. Indeed, as we have discussed before, Blasco \cite[Theorem 1.2]{Blasco-OTAA-2010} has also observed that $K(\ell^2_q)=0$ for every $1\leq q\leq \infty$.
\end{rem}

Our main aim of this paper is to investigate precise asymptotic behavior of $\lambda_{p}$-Bohr radius $K^p(B_{\ell^n_q},T,\lambda)$, where $T$ is a bounded operator between two complex Banach spaces. Motivated by the work of Defant {\it et al.} \cite{defant-adv-math-2012}, one of our contributions in this paper is the following theorem, which provides the correct asymptotic estimates for $\lambda_{p}$-Bohr radius of identity operator on Banach spaces, where $1\leq p<\infty$. Let us first concentrate on finite dimensional Banach spaces. 
\begin{thm}\label{Pal-Vasu-P5-thm-01}
	Let $X$ be a finite dimensional complex Banach space and $\lambda>1$, $1\leq p<\infty$. Then there are constants $B$ and $C$ such that for each $1\leq q\leq \infty$ and $n\in \mathbb{N}$ we have 
	\begin{equation*}
		K^p(B_{\ell^n_q},X,\lambda)\leq B \lambda^{2/p}n^{1-1/p}\left(\frac{\log n}{n}\right)^{1-\frac{1}{\min\{q, 2\}}}
	\end{equation*}
	and
	$$
	K^p(B_{\ell^n_q},X,\lambda)\geq \left\{\begin{array}{ll}
		C\left(\frac{\lambda-1}{2\lambda-1}\right)^{1/p} \left(\frac{\log n}{n}\right)^{1-\frac{1}{r}}, & \mbox{for \, $q\leq r$}\\[5mm] 
		
		C\left(\frac{\lambda-1}{2\lambda-1}\right)^{1/p} n^{\frac{1}{q}-\frac{1}{r}}\left(\frac{\log n}{n}\right)^{1-\frac{1}{r}}, & \mbox{for \, $r<q$}
	\end{array}\right.
	$$
	where $r=pq'/(pq'-1)$. Here $B$ is a universal constant and $C=e^{-1}\inf_{m\in \mathbb{N}}(\pi_{1}(\mbox{Id}_{X}))^{-1/m}$, depends only on $X$.
\end{thm}
Now, we want to study the Bohr radius $K^p(B_{\ell^n_q},X,\lambda)$ for the infinite dimensional Banach spaces. In this cases, we observe that the log term is missing from the estimates of $K^p(B_{\ell^n_q},X,\lambda)$. In fact, the asymptotic decay is very much influenced by geometry of the Banach space $X$, {\it i.e.,} by its optimal cotype $\text{Cot}(X)$.

\begin{thm}\label{Pal-Vasu-P5-thm-02}
	Let $1\leq p<\infty$ and $1\leq q\leq \infty$. If $X$ is an infinite dimensional complex Banach space of cotype $r$, then we have 
	$$
	K^p(B_{\ell^n_q},X,\lambda)\leq \left\{\begin{array}{ll}
		\frac{\lambda^{1/p}}{n^{\frac{1}{p}-\frac{1}{q}}}, & \mbox{for \, $q\leq \text{Cot}(X)$}\\[5mm] 
		
		\frac{\lambda^{1/p}}{n^{\frac{1}{p}-\frac{1}{\text{Cot}(X)}}}, & \mbox{for \, $q>\text{Cot}(X)$}
	\end{array}\right.
	$$
	and 
	$$
	K^p(B_{\ell^n_q},X,\lambda)\geq \left\{\begin{array}{ll}
		\left(\frac{\lambda-1}{\lambda}\right)^{1/p}\frac{1}{e\, C_r(X)}, & \mbox{for \, $q\leq s$}\\[5mm] 
		
		\left(\frac{\lambda-1}{\lambda}\right)^{1/p}\frac{1}{e\,C_r(X)n^{\frac{1}{s}-\frac{1}{q}}}, & \mbox{for \, $q>s$},
	\end{array}\right.
	$$
	where $s=r/(r-1)$.\\
	
	\noindent
	Moreover, if $X$ has no finite cotype, then 
	\begin{equation*}
		\left(\frac{\lambda-1}{\lambda}\right)^{\frac{1}{p}}\frac{1}{n^{\frac{1}{p}-\frac{1}{pq}}} \leq K^p(B_{\ell^n_q},X,\lambda) \leq \frac{\lambda^{1/p}}{n^{\frac{1}{p}-\frac{1}{q}}}.
	\end{equation*}
\end{thm}
\vspace{4mm}

Defant {\it et al.} \cite[Theorem 4.1]{defant-adv-math-2012} have used the hypercontractivity of the Bohnenblust-Hille inequality to establish the asymptotic behavior of of the Bohr radii $K^p(B_{\ell^n_{\infty}},X, \lambda)$ for $p=1$. Moreover, as a technical tool, Defant {\it et al.} \cite[Theorem 5.3]{defant-adv-math-2012} have shown Bohnenblust-Hille type inequality for $(s,1)$-summing operator to estimate the multidimensional Bohr radii for operators: {\it let $Y$ be a $r$-concave Banach lattice, with $2\leq r<\infty$, and $T:X\rightarrow Y$ an $(s,1)$-summing operator with $1\leq s\leq r$. Then there is a constant $C>0$ such that for every $m$-homogeneous polynomial $Q:\ell^n_{\infty}\rightarrow X$, $Q(z)=\sum_{\alpha \in\Lambda(m,n)}x_{\alpha}z^{\alpha}$}

\begin{equation*}
	\left(\sum_{\alpha \in\Lambda(m,n)}\norm{T(x_{\alpha})}^{\rho}_{Y}\right)^{\frac{1}{\rho}}\leq C^m\norm{Q}_{B_{\ell^n_{\infty}}},
\end{equation*}
where 
\begin{equation*}
	\rho:=\frac{rsm}{r+(m-1)s}.
\end{equation*}
and
\begin{equation*}
	C^m=\pi_{s,1}(T)M_{r}(Y)(\sqrt{2})^{m-1}m^{\frac{r-1}{r}}\left(\frac{m}{m-1}\right)^{m-1}.
\end{equation*}
Recently, by considering the diagonal operator $D_{z}:\ell^n_{\infty}\rightarrow \ell^n_{q}$, Kumar and Manna \cite{kumar-2023-arxiv} have obtained an extended version of the above result for $m$-homogeneous polynomials $Q\in \mathcal{P}(^m\ell^n_q,X), \,\, 1\leq q<\infty$.

\begin{thm}\cite{kumar-2023-arxiv}\label{Kumar-thm-arxiv}
	Let $2\leq r<\infty$. Suppose $Y$ is a $r$-concave Banach lattice, and the operator $T:X\rightarrow Y$ is an $(s,1)$-summing operator with $1\leq s\leq r$. Then there is a constant $\widetilde{C}>0$ such that the following holds
	\begin{equation*}
		\left(\sum_{\alpha \in\Lambda(m,n)}\norm{T(x_{\alpha})z^{\alpha}}^{\rho}_{Y}\right)^{\frac{1}{\rho}}\leq \widetilde{C}^m\norm{Q}_{B_{\ell^n_{q}}}
	\end{equation*}
	for every polynomial $Q(z)=\sum_{\alpha \in\Lambda(m,n)}x_{\alpha}z^{\alpha}\in \mathcal{P}(^m\ell^n_{q}),\,\, 1\leq q<\infty$.
\end{thm}

We are now ready to state our next main result. As a consequence of the above technical discussions, Theorem \ref{Pal-Vasu-P5-thm-03} provides the lower estimates of the Bohr radius $K^p(B_{\ell^n_q},T,\lambda)$ for any bounded operator $T$.
\begin{thm}\label{Pal-Vasu-P5-thm-03}
	Let $1\leq p<\infty$ and $T:X\rightarrow Y$ be a bounded operator between Banach spaces $X$ and $Y$.
	\begin{enumerate}
		\item Assume that $X$ or $Y$ is of cotype $r$ with $2\leq r\leq \infty$. Then there is a constant $D>0$ such that for every $\norm{T}<\lambda^{1/p}$, $1\leq q\leq \infty$, and $n$
		$$
		K^p(B_{\ell^n_q},T, \lambda) \geq \left\{\begin{array}{ll}
			D\, \left(\frac{\lambda - \norm{T}^p}{\lambda}\right)^{\frac{1}{p}}, & \mbox{for \, $q\leq s$}\\[5mm] 
			
			D\, \left(\frac{\lambda - \norm{T}^p}{\lambda}\right)^{\frac{1}{p}} \left(\frac{1}{n}\right)^{\frac{1}{s}-\frac{1}{q}}, & \mbox{for \, $q>s$},
		\end{array}\right.
		$$
		where $s=r/(r-1)$ and $D=(e\, \min\{C_r(X),C_r(Y)\})^{-1}$.\\[2mm]
		\item Assume that $Y$ is a $r$-concave Banach lattice with $2\leq r<\infty$ and there is a $1\leq s<r$ such that the operator $T$ is $(s,1)$-summing. Then there is a constant $D>0$  such that 
		\begin{equation*}
			K^p(B_{\ell^n_q},T, \lambda) \geq D\, \left(\frac{\lambda-\norm{T}^p}{2\lambda -\norm{T}^p}\right)^{\frac{1}{p}} \left(\frac{\log n}{n}\right)^{1-\frac{1}{r}}
		\end{equation*}
	for every $\norm{T}<\lambda^{1/p}$, $1\leq q\leq \infty$, and $n\in \mathbb{N}$ and the constant $D$ is given by
	\begin{equation*}
		D=\left(\frac{1}{2}\right)^{1-\frac{1}{r}}e \left(\pi_{s,1}(T)M_r(Y)\right)^{1/m}(\sqrt{2})^{1-\frac{1}{m}}m^{\frac{r-1}{rm}}\left(\frac{m}{m-1}\right)^{1-\frac{1}{m}}.
	\end{equation*}
	\end{enumerate}
\end{thm}
We now move our focus for a while on arithmetic Bohr radius of Banach spaces, which is another interesting multidimensional variant of Bohr radius. In fact, it has very close relationship with classical multidimensional Bohr radii. In order to study the upper inclusion of domains of convergence of monomial expansions, as a technical tool, Defant {\it et al.} \cite{defant-Domain-J angew-2009} introduced the arithmetic Bohr radius for scalar-valued holomorphic functions. Later, it was studied independently by the same authors in \cite{defant-QJM-2008}. Recently, the study of arithmetic Bohr radius has been extended for vector-valued holomorphic functions (see \cite{Subhadip-Vasu-P3}). Indeed, the authors \cite{Subhadip-Vasu-P3} have introduced the following notion of arithmetic Bohr radius:
\begin{defn}\cite{Subhadip-Vasu-P3}
	Let $\mathcal{F}(B_{\ell^n_q},X)$ be set of all holomorphic functions on $B_{\ell^n_q}$ into a complex Banach space $X$. For each $1\leq p<\infty$ and $\lambda\geq 1$, the $\lambda_{p}$-\textit{arithmetic Bohr radius} of $B_{\ell^n_q}$ with respect to $\mathcal{F}(B_{\ell^n_q}, X)$ is defined as 
	\begin{equation*}
		A_{p}(\mathcal{F}(B_{\ell^n_q}, X), \lambda) := \sup \left\{\frac{1}{n}\sum_{i=1}^{n}r_i \,|\, r\in \mathbb{R}^{n}_{\geq 0},\, \forall\, f\in \mathcal{F}(B_{\ell^n_q}, X) : \sum_{\alpha \in \mathbb{N}^{n}_{0}}\norm{x_{\alpha}(f)}^p r^{p\alpha} \leq \lambda \norm{f}^{p}_{B_{\ell^n_q}}\right\},
	\end{equation*} 
	where $\mathbb{R}^n_{\geq 0}=\{r=(r_1,.\,.\,.\,, r_n) \in \mathbb{R}^n: r_i\geq 0, 1\leq i\leq n\}.$  We write $A_p(B_{\ell^n_q},X,\lambda)$ for $A_p(H_{\infty}(B_{\ell^n_q},X),\lambda)$ and $A_p(B_{\ell^n_q},X)$ for $A_p(B_{\ell^n_q},X,1)$.
\end{defn}
In view of \cite[Theorem 2.2]{Subhadip-Vasu-P3} and Theorem \ref{Pal-Vasu-P5-thm-01}, as an application, we obtain the following lower estimate for the arithmetic Bohr radius $A_p(B_{\ell^n_q},X,\lambda)$ in case of $X$ is finite dimensional.
\begin{thm}
	Let $1\leq p<\infty$, $\lambda\geq 1$, and $1\leq q\leq \infty$. Suppose $X$ is a finite dimensional Banach space. Then there is a constant $C>0$ such that for every $n>1$ we have
	$$
	A_{p}(B_{\ell^n_q},X,\lambda)\geq \left\{\begin{array}{ll}
		C\left(\frac{\lambda-1}{2\lambda-1}\right)^{1/p} n^{-\frac{1}{q}}\left(\frac{\log n}{n}\right)^{1-\frac{1}{r}}, & \mbox{for \, $q\leq r$}\\[5mm] 
		
		C\left(\frac{\lambda-1}{2\lambda-1}\right)^{1/p} n^{-\frac{1}{r}}\left(\frac{\log n}{n}\right)^{1-\frac{1}{r}}, & \mbox{for \, $r<q$},
	\end{array}\right.
	$$
	where $r=pq'/(pq'-1)$, and the constant $C=e^{-1}\inf_{m\in \mathbb{N}}(\pi_{1}(\mbox{Id}_{X}))^{-1/m}$, depends only on $X$.
\end{thm}
\section{Key Lemmas}\label{Pal-Vasu-P5-section-04}
\allowdisplaybreaks
As demonstrated by the scalar case, studying the $m$-homogeneous polynomial first is typically a beneficial approach for obtaining nontrivial estimates for Bohr radii. The following definition is the $m$-homogeneous analogue of Definition \ref{Pal-Vasu-P5-defn-01} for $m\in \mathbb{N}$.

\begin{defn}\label{Pal-Vasu-P5-defn-02}
	Let $1\leq p<\infty$, $1\leq q\leq \infty$, $\lambda \geq 1$, $m,n \in \mathbb{N}$, and $T: X\rightarrow Y$ be a bounded operator in complex Banach spaces such that $\norm{T}\leq \lambda^{1/p}$. Then $K^p_m(B_{\ell^n_q},T,\lambda)$ is defined to be the supremum of all $r\geq 0$ such that for every $m$-homogeneous polynomials $Q\in \mathcal{P}(^m\ell^n_q, X)$, $Q(z)=\sum_{\alpha \in\Lambda(m,n)}x_{\alpha}z^{\alpha}$, we have 
	\begin{equation*}
		\sup_{z\in rB_{\ell^n_{q}}}\sum_{\alpha \in\Lambda(m,n)}\norm{T(x_{\alpha})z^{\alpha}}^{p}_{Y}\leq \lambda\sup_{z\in B_{\ell^n_{q}}}\norm{\sum_{\alpha \in\Lambda(m,n)}x_{\alpha}z^{\alpha}}^p_{X}.
	\end{equation*}
\end{defn}
It is worth mentioning that if $T$ is a null operator, then the above supremum is $\infty$. Thus, by assuming that $T\neq 0$ throughout the paper, we are avoiding this trivial scenario. The constant $K^p_m(B_{\ell^n_q},T,\lambda)$ in the Definition \ref{Pal-Vasu-P5-defn-02} can be rephrased as
\begin{align}\label{Pal-Vasu-P5-eqn-int-01}
	K^p_m(B_{\ell^n_q},T,\lambda):= \sup \biggl\{ & r\geq 0 : \sup_{z\in B_{\ell^n_{q}}}\sum_{\alpha \in\Lambda(m,n)}\norm{T(x_{\alpha})z^{\alpha}}^{p}_{Y}\leq \frac{\lambda}{r^{mp}}\sup_{z\in B_{\ell^n_{q}}}\norm{Q(z)}^p_{X}\\[2mm]\nonumber & \mbox{for all} \,\, Q\in \mathcal{P}(^m\ell^n_q, X)\biggr\}.
\end{align}
Further, it is worth noting that
\begin{equation}\label{Pal-Vasu-P5-eqn-int-02}
	K^p_m(B_{\ell^n_q},T,\lambda)=\sqrt[mp]{\lambda}K^p_m(B_{\ell^n_q},T)
\end{equation}
and 
\begin{equation}\label{Pal-Vasu-P5-eqn-int-03}
	K^p(B_{\ell^n_q},T,\lambda) \geq \max \left\{K^p\left(B_{\ell^n_q},X,\lambda/\norm{T}^p\right), K^p\left(B_{\ell^n_q},Y,\lambda/\norm{T}^p\right)\right\}.
\end{equation}
The following lemma provides a nice connection between Definitions \ref{Pal-Vasu-P5-defn-01} and \ref{Pal-Vasu-P5-defn-02}, which will be one of the most important technical tools to prove our main results in this paper.
\begin{lem}\label{Pal-Vasu-P5-Lemma-01}
	Let  $1\leq p<\infty$ and $T:X\rightarrow Y$ be a bounded operator between complex Banach spaces $X$ and $Y$ with $\norm{T}<\lambda^{1/p}$. Then we have 
	\begin{enumerate}
		\item[(a)] $\left(\frac{\lambda-\norm{T}^p}{2\lambda-\norm{T}^p}\right)^{1/p}\inf_{m\in \mathbb{N}}K^p_m(B_{\ell^n_q}, T, \lambda)\leq K^p(B_{\ell^n_q},T,\lambda)\leq \inf_{m\in \mathbb{N}}K^p_m(B_{\ell^n_q},T,\lambda)$
		\item[(b)] $\left(\frac{\lambda-\norm{T}^p}{\lambda-\norm{T}^p+1}\right)^{1/p}\inf_{m\in\mathbb{N}}K^p_m(B_{\ell^n_q},T)\leq K^p(B_{\ell^n_q},T,\lambda)\leq \lambda^{1/p} \inf_{m\in \mathbb{N}}K^p_m(B_{\ell^n_q},T).$
	\end{enumerate}
\end{lem}

\begin{pf}
	The right hand side inequality of (a) is clear from the fact that $K^p(B_{\ell^n_{q}},T, \lambda)\leq K^p_{m}(B_{\ell^n_{q}},T, \lambda)$ for all $m\geq 1$. So, we now proceed for the left hand side inequality of (a). Let $f(z)=\sum_{\alpha\in \mathbb{N}^n_{0}}x_{\alpha}z^{\alpha}$ be a holomorphic function on $B_{\ell^n_q}$ with values in the complex Banach space $X$. Further, suppose $v\in \mathbb{C}^n$ such that $\norm{v}\leq \inf_{m\in \mathbb{N}}\{K^p_{m}(B_{\ell^n_{q}},T, \lambda)\}.$
	It is worth noting that for $1\leq p<\infty$
	\begin{equation*}
		0<	\left(\frac{\lambda-\norm{T}^p}{2\lambda-\norm{T}^p}\right)^{1/p}<1.
	\end{equation*}
	Since $v\in K^p_{m}(B_{\ell^n_{q}},T,\lambda)\overline{B_{\ell^n_{q}}}$ for all $m\in \mathbb{N}$ and using the Cauchy-Riemann inequalities we have that
	\begin{align*}
		&\sum_{\alpha\in \mathbb{N}^n_{0}}\norm{T(x_{\alpha})\left(\left(\frac{\lambda-\norm{T}^p}{2\lambda-\norm{T}^p}\right)^{1/p} v\right)^{\alpha}}^p_{Y} \\[2mm]&= \norm{T(x_0)}^p+\sum_{m=1}^{\infty}\left(\frac{\lambda-\norm{T}^p}{2\lambda-\norm{T}^p}\right)^{m}\sum_{\alpha \in\Lambda(m,n)}\norm{T(x_{\alpha})}^p|v|^{p\alpha} \\[2mm]&
		\leq \norm{T}^p\norm{x_0}^p + \sum_{m=1}^{\infty}\left(\frac{\lambda-\norm{T}^p}{2\lambda-\norm{T}^p}\right)^{m} \lambda \norm{\sum_{\alpha \in\Lambda(m,n)}x_{\alpha}z^{\alpha}}^p_{B_{\ell^n_q}}\\[2mm]&
		\leq \left(\norm{T}^p + \lambda \sum_{m=1}^{\infty} \left(\frac{\lambda-\norm{T}^p}{2\lambda-\norm{T}^p}\right)^m\right)\norm{f}^p_{B_{\ell^n_q}}=\lambda \norm{f}^p_{B_{\ell^n_q}}.
	\end{align*}
	We shall prove (b) in the similar way. Since $K^p_m(B_{\ell^n_q},T,\lambda)=\lambda^{1/pm}K^p_m(B_{\ell^n_q},T)$, we obtain 
	\begin{equation*}
		K^p(B_{\ell^n_q},T,\lambda)\leq \inf_{m\in \mathbb{N}} \left\{\lambda^{\frac{1}{pm}}K^p_m(B_{\ell^n_q},T)\right\}\leq \lambda^{1/p} \inf_{m\in \mathbb{N}}\{K^p_m(B_{\ell^n_q},T)\}.
	\end{equation*}
	In this case, we suppose again that $v\in \mathbb{C}^n$ as above but here we proceed by considering $\norm{v}\leq \inf_{m\in \mathbb{N}} \{K^p_m(B_{\ell^n_q},T)\}$. Thus, for $v\in K^p_m(B_{\ell^n_q},T)\overline{B_{\ell^n_{q}}}$, we obtain
	\begin{align*}
		&	\sum_{\alpha\in \mathbb{N}^n_{0}} \norm{T(x_{\alpha})\left(\left(\frac{\lambda-\norm{T}^p}{\lambda-\norm{T}^p+1}\right)^{1/p} v\right)^{\alpha}}^p_{Y}\\[2mm]&=
		\norm{T(x_0)}^p+\sum_{m=1}^{\infty}\left(\frac{\lambda-\norm{T}^p}{\lambda-\norm{T}^p+1}\right)^{m}\sum_{\alpha \in\Lambda(m,n)}\norm{T(x_{\alpha})}^p|v|^{p\alpha}\\[2mm]&
		\leq \norm{T}^p\norm{x_0}^p + \sum_{m=1}^{\infty}\left(\frac{\lambda-\norm{T}^p}{\lambda-\norm{T}^p+1}\right)^{m}  \norm{\sum_{\alpha \in\Lambda(m,n)}x_{\alpha}z^{\alpha}}^p_{B_{\ell^n_q}}\\[2mm]&
		\leq \left(\norm{T}^p +  \sum_{m=1}^{\infty} \left(\frac{\lambda-\norm{T}^p}{\lambda-\norm{T}^p+1}\right)^m\right)\norm{f}^p_{B_{\ell^n_q}}=\lambda \norm{f}^p_{B_{\ell^n_q}}.
	\end{align*}
	This completes the proof.
\end{pf}

The following lemma of Bayart {\it et al.} \cite[Lemma 3.5]{bayart-J-Anal-2019} is a key tool to prove Theorem \ref{Pal-Vasu-P5-thm-01}. For a $m$-homogeneous polynomial $Q\in \mathcal{P}(^m\ell^n_q), \,\,(1\leq q\leq \infty)$, it guarantees the $\ell_{q'}$-summability of certain coefficients of $Q$ in terms of its uniform norm on $\ell^n_q$.

\begin{lem}\cite{bayart-J-Anal-2019}\label{Pal-Vasu-P5-Lemma-02}
	Let $1\leq t\leq \infty$ and $P$ be an $m$-homogeneous polynomial in $n$ variables. Then for any $\mathbf{j}\in \mathcal{J}(m-1, n)$,
	\begin{equation*}
		\left(\sum_{k=j_{m-1}}^{n}|c_{(\mathbf{j},k)}(P)|^{t'}\right)^{1/t'}\leq me^{1+\frac{m-1}{t}}|\mathbf{j}|^{1/t}\norm{P}_{\mathcal{P}(^m\ell^n_t)}.
	\end{equation*}
\end{lem}

\section{Proof of the main results}\label{Pal-Vasu-P5-section-05}

\begin{pf}[{\bf Proof of Proposition \ref{Pal-Vasu-P5-prop-01}}]
	Let $Q(z)=\sum_{\alpha \in\Lambda(m,n)}x_{\alpha}z^{\alpha}\in \mathcal{P}(^ml^n_q,X)$ be an $m$-homogeneous polynomial. By using H\"{o}lder's inequality and for all $z\in \mathbb{C}^n$, we have 
	\begin{align*}
		\sum_{\alpha \in\Lambda(m,n)}\norm{T(x_{\alpha})z^{\alpha}}^p &\leq \max_{\alpha \in\Lambda(m,n)}\norm{T(x_{\alpha})}^p \sum_{\alpha \in\Lambda(m,n)}|z^{\alpha}|^p \\[2mm]& \leq \norm{T}^p \max_{\alpha \in\Lambda(m,n)}\norm{x_{\alpha}}^p \left(|z_1|^p+\cdots+|z_n|^p\right)^m \\[2mm]& 
		\leq \lambda \max_{\alpha \in\Lambda(m,n)}\norm{x_{\alpha}}^p \left(|z_1|^{pq}+\cdots+|z_n|^{pq}\right)^{\frac{m}{q}} n^{(1-\frac{1}{q})m}\\[2mm]& 
		\leq \lambda \norm{Q}^p_{\ell^n_q} \norm{z}^{mp}_{pq}n^{{mp}\left(\frac{1}{p}-\frac{1}{pq}\right)}\\[2mm]&
		\leq \lambda \norm{Q}^p_{\ell^n_q} \norm{z}^{mp}_{q}n^{{mp}\left(\frac{1}{p}-\frac{1}{pq}\right)}.
	\end{align*}
Therefore, for all $z\in \left(1/n^{(1/p-1/pq)}\right)B_{\ell^n_q}$, we obtain
\begin{equation*}
	\sum_{\alpha \in\Lambda(m,n)}\norm{T(x_{\alpha})z^{\alpha}}^p\leq \lambda \norm{Q}^p.
\end{equation*}
Thus, for every $m$-homogeneous polynomial $Q\in \mathcal{P}(^m\ell^n_q)$, we obtain
\begin{equation*}
	\sup \left\{\sum_{\alpha \in\Lambda(m,n)}\norm{T(x_{\alpha})z^{\alpha}}^p \,| \,z\in \frac{1}{n^{\frac{1}{p}-\frac{1}{pq}}}B_{\ell^n_q}\right\}\leq \lambda \norm{Q}^p,
\end{equation*}
which yields 
\begin{equation*}
	K^p_m(B_{\ell^n_q},T,\lambda)\geq \frac{1}{n^{\frac{1}{p}-\frac{1}{pq}}}.
\end{equation*}
In view of Lemma \ref{Pal-Vasu-P5-Lemma-01}(a), for each $1\leq p<\infty$, we have 
\begin{equation*}
	K^p(B_{\ell^n_q},T,\lambda)\geq \frac{1}{n^{\frac{1}{p}-\frac{1}{pq}}} \left(\frac{\lambda-\norm{T}^p}{2\lambda-\norm{T}^p}\right)^{1/p}.
\end{equation*}
By going similar lines of arguments as above for every $m$-homogeneous polynomial $Q(z)=\sum_{|\alpha|=m}x_{\alpha}z^{\alpha}$, we obtain that
\begin{equation*}
	\sup_{z\in \left(\frac{1}{\sqrt[m]{\norm{T}}n^{\frac{1}{p}-\frac{1}{pq}}}\right)B_{\ell^n_q}} \sum_{\alpha \in\Lambda(m,n)}\norm{T(x_{\alpha})z^{\alpha}}^p=
	\frac{1}{\norm{T}^p} \sup_{z\in \left(\frac{1}{n^{\frac{1}{p}-\frac{1}{pq}}}\right)B_{\ell^n_q}} \sum_{\alpha \in\Lambda(m,n)}\norm{T(x_{\alpha})z^{\alpha}}^p\leq \norm{Q}^p_{\mathcal{P}(^m\ell^n_q)}.
\end{equation*}
Hence, for every $m\in \mathbb{N}$, we have 
\begin{equation*}
	K^p_m(B_{\ell^n_q},T)\geq \frac{1}{\sqrt[m]{\norm{T}}n^{\frac{1}{p}-\frac{1}{pq}}}.
\end{equation*}
Now, by virtue of Lemma \ref{Pal-Vasu-P5-Lemma-01}(b), we have 
\begin{equation*}
	K^p(B_{\ell^n_q},T,\lambda)\geq \left(\frac{\lambda-\norm{T}^p}{\lambda-\norm{T}^p+1}\right)^{1/p} \frac{1}{\norm{T}n^{\frac{1}{p}-\frac{1}{pq}}} \quad \mbox{if}\,\, \norm{T}\geq 1
\end{equation*}
and 
\begin{equation*}
	K^p(B_{\ell^n_q},T,\lambda)\geq \left(\frac{\lambda-\norm{T}^p}{\lambda-\norm{T}^p+1}\right)^{1/p} \frac{1}{n^{\frac{1}{p}-\frac{1}{pq}}} \quad \mbox{if}\,\, 0<\norm{T}< 1.
\end{equation*}
This completes the proof.
\end{pf}

\begin{pf}[{\bf Proof of Theorem \ref{Pal-Vasu-P5-thm-01}}]
	We first check that the upper bound of $K^p(B_{\ell^n_q},X,\lambda)$ is true. In view of \cite[Lemma 2.3]{Subhadip-Vasu-P3}, we have 
	\begin{equation}\label{Pal-Vasu-P5-e-001}
		K^p(B_{\ell^n_q},X,\lambda)\leq n^{1/q}A_p(B_{\ell^n_q},X,\lambda),\,\, \lambda\geq 1
	\end{equation}
for any $1\leq p<\infty$ and $n\in \mathbb{N}$, where $A_p(B_{\ell^n_q},X,\lambda)$ is the $\lambda_p$-arithmetic Bohr radius for the unit ball of $\ell^n_q$ spaces.
In the same paper \cite[Theorem 2.2]{Subhadip-Vasu-P3}, the authors have established that there exists a universal constant $B>0$ such that 
\begin{equation}\label{Pal-Vasu-P5-e-002}
	 A_p(B_{\ell^n_q},X, \lambda)\leq B \,\lambda^{\frac{2}{p\log n}}\, \frac{\left(\log n\right)^{1-\left(1/\min\{q, 2\}\right)}}{n^{\frac{1}{p}+\frac{1}{\max\{2, q\}}-\frac{1}{2} }}.
\end{equation}
Using \eqref{Pal-Vasu-P5-e-001} and \eqref{Pal-Vasu-P5-e-002} and a simple computation yields 
\begin{equation*}
	K^p(B_{\ell^n_q},X,\lambda) \leq B \lambda^{\frac{2}{p\log n}}\,\,\frac{(\log n)^{1-\frac{1}{\min \{q,2\}}}}{n^{\frac{1}{p}-\frac{1}{\min\{q,2\}}}}=B\,\lambda^{\frac{2}{p\log n}}\, n^{1-\frac{1}{p}} \left(\frac{\log n}{n}\right)^{1-\frac{1}{\min\{q,2\}}}.
\end{equation*}
\allowdisplaybreaks
Now we shall prove the lower bound for $K^p(B_{\ell^n_q},X,\lambda)$. By virtue of Lemma \ref{Pal-Vasu-P5-Lemma-01}, it is sufficient to show the existence of a constant $C(X)>0$ such that 
	$$
\inf_{m\in \mathbb{N}} K^p_m(B_{\ell^n_q},X,\lambda)\geq \left\{\begin{array}{ll}
	C \,\left(\frac{\log n}{n}\right)^{1-\frac{1}{r}}, & \mbox{for \, $q\leq r$}\\[6mm] 
	
	C \,n^{\frac{1}{q}-\frac{1}{r}}\left(\frac{\log n}{n}\right)^{1-\frac{1}{r}}, & \mbox{for \, $r<q$}
\end{array}\right.
$$
for every $1\leq p<\infty$.
Let us fix $Q\in \mathcal{P}(^m\ell^n_q,X)$ and $u\in \ell^n_q$. Then, by Lemma \ref{Pal-Vasu-P5-Lemma-02}, it is worth noting that
\begin{align}\label{Pal-Vasu-P5-e-003}
	\left(\sum_{k: \,\,(\mathbf{j},k)\in \mathcal{J}(m,n)}|c_{(\mathbf{j},k)}(P)|^{t'}\right)^{1/t'}\nonumber &=\left(\sum_{k=j_{m-1}}^{n}|c_{(\mathbf{j},k)}(P)|^{t'}\right)^{1/t'}\\[3mm] & 
	\leq me^{1+\frac{m-1}{t}}|\mathbf{j}|^{1/t}\norm{P}_{\mathcal{P}(^m\ell^n_t)} 
\end{align}
holds for any $m$-homogeneous polynomial $P\in \mathcal{P}(^m\ell^n_t), \, 1\leq t\leq \infty$.
Since $X$ is a finite dimensional Banach space, the identity map $\mbox{Id}_X$ on $X$ is absolutely summing {\it i.e.}, $1$-summing and hence the operator $\mbox{Id}_X$ is $s$-summing for every $s\geq 1$ with $\pi_s(\mbox{Id}_X)\leq \pi_1(\mbox{Id}_X)$ (see \cite[Theorem 10.4]{diestel-Abs-sum-book}). In other words, there exists a constant $d>0$ such that for every $s\geq 1$ and for every finitely many $x_1, \ldots, x_k \in X$ ,
\begin{equation*}
	\left(\sum_{i=1}^{k}\norm{x_i}^s\right)^{1/s}\leq d \sup_{x^*\in B_{X^*}}\left(\sum_{i=1}^{k}|x^*(x_i)|^s\right)^{1/s}.
\end{equation*}
Hence, if $Q(z)=\sum_{|\alpha|=m}x_{\alpha}(Q)z^{\alpha}$, then using the above inequality \eqref{Pal-Vasu-P5-e-003}, H\"{o}lder's inequality (two times), and the multinomial formula, for every $1\leq p<\infty$ we obtain
\begin{align*}
&\sum_{\mathbf{j}\in \mathcal{J}(m,n)} \norm{x_{\mathbf{j}}(Q)}^p|u_{\mathbf{j}}|^p= \sum_{\mathbf{j}\in \mathcal{J}(m,n)^*}\left(\sum_{k: \,\,(\mathbf{j},k)\in \mathcal{J}(m,n)}\norm{x_{\mathbf{j}}(Q)}^p|u_{\mathbf{j}}|^p|u_k|^p\right)\\[3mm] & \leq \sum_{\mathbf{j}\in \mathcal{J}(m,n)^*}|u_{\mathbf{j}}|^p \left(\sum_{k: \,\,(\mathbf{j},k)\in \mathcal{J}(m,n)}\norm{x_{\mathbf{j}}(Q)}^{pq'}\right)^{1/q'} \left(\sum_{k}|u_k|^{pq}\right)^{1/q}\\[3mm] &
=\sum_{\mathbf{j}\in \mathcal{J}(m,n)^*}|u_{\mathbf{j}}|^p \left(\left(\sum_{k: \,\,(\mathbf{j},k)\in \mathcal{J}(m,n)}\norm{x_{\mathbf{j}}(Q)}^{pq'}\right)^{1/pq'}\norm{u}_{pq}\right)^p\\[3mm] & \leq \pi^p_1(\mbox{Id}_X) \sum_{\mathbf{j}\in \mathcal{J}(m,n)^*}|u_{\mathbf{j}}|^p \left(\sup_{x^*\in B_{X^*}}\left(\sum_{k: \,\,(\mathbf{j},k)\in \mathcal{J}(m,n)}|x^*(x_{(\mathbf{j},k)}(Q))|^{pq'}\right)^{1/pq'}\right)^p \norm{u}^p_{pq}\\[3mm] & 
\leq \pi^p_1(\mbox{Id}_X) \left(me^{1+\frac{m-1}{r}}\right)^p \norm{u}^p_{pq} \sup_{x^*\in B_{X^*}}\sup_{z\in \ell^n_q} \bigg|x^*\left(\sum_{\mathbf{j}\in \mathcal{J}(m,n)}x_{\alpha}(Q)z^{\alpha}\right)\bigg|^p \sum_{\mathbf{j}\in \mathcal{J}(m,n)^*} |\mathbf{j}|^{p/r}|u_{\mathbf{j}}|^p\\[3mm]&
\leq \pi^p_1(\mbox{Id}_X) \left(me^{1+\frac{m-1}{r}}\right)^p \norm{u}^p_q \norm{Q}^p_{\mathcal{P}(^m\ell^n_q)}\left(\sum_{\mathbf{j}\in \mathcal{J}(m,n)^*}|\mathbf{j}|^{1/r}|u_{\mathbf{j}}|\right)^p\\[3mm]&
\leq \pi^p_1(\mbox{Id}_X) \left(me^{1+\frac{m-1}{r}}\right)^p \norm{u}^p_q \norm{Q}^p_{\mathcal{P}(^m\ell^n_q)} \left(\left(\sum_{\mathbf{j}\in \mathcal{J}(m,n)^*}|\mathbf{j}||u_{\mathbf{j}}|^r\right)^{1/r}\left(\sum_{\mathbf{j}\in \mathcal{J}(m,n)^*}1\right)^{1-1/r}\right)^p.
\end{align*}
It is worth noting that 
\begin{equation*}
	\left(\sum_{\mathbf{j}\in \mathcal{J}(m,n)^*}|\mathbf{j}||u_{\mathbf{j}}|^r\right)^{1/r}=\norm{u}^{m-1}_r.
\end{equation*}
Hence, we obtain 
\begin{equation*}
	\sum_{\mathbf{j}\in \mathcal{J}(m,n)} \norm{x_{\mathbf{j}}(Q)}^p|u_{\mathbf{j}}|^p \leq \pi^p_1(\mbox{Id}_X) \left(me^{1+\frac{m-1}{r}}\right)^p \norm{u}^{p(m-1)}_r\norm{u}^p_{q}\norm{Q}^p_{\mathcal{P}(^m\ell^n_q)}\left(\sum_{\mathbf{j}\in \mathcal{J}(m,n)^*}1\right)^{p\left(1-\frac{1}{r}\right)}.
\end{equation*}
Using the facts that 
\begin{equation*}
	\mbox{card} \,\{\mathbf{j}=(j_1,\ldots,j_{m-1}): 1\leq j_1\leq \cdots \leq j_{m-1}\leq n\}= \binom{n+m-2}{n-1}
\end{equation*}
and 
\begin{equation*}
	\binom{n+m-2}{n-1}\leq e^m\left(1+\frac{n}{m}\right)^{m-1} \,\,\, \mbox{for all}\,\, n,m\in \mathbb{N},
\end{equation*}
finally we have the following inequality
\begin{equation*}
	\sum_{\mathbf{j}\in \mathcal{J}(m,n)} \norm{x_{\mathbf{j}}(Q)}^p|u_{\mathbf{j}}|^p \leq \pi^p_1(\mbox{Id}_X) \left(m e^{m+\frac{1}{r'}}\right)^p \left(1+\frac{n}{m}\right)^{\frac{p(m-1)}{r'}} \norm{u}^{p(m-1)}_{r}\norm{u}^p_{q} \norm{Q}^p_{\mathcal{P}(^m\ell^n_q)}.
\end{equation*}
For $q\leq r$, we have $\norm{u}_r\leq \norm{u}_q$. Then, we obtain
\begin{equation}\label{Pal-Vasu-P5-eqn-004}
		\sum_{\mathbf{j}\in \mathcal{J}(m,n)} \norm{x_{\mathbf{j}}(Q)}^p|u_{\mathbf{j}}|^p \leq \pi^p_1(\mbox{Id}_X) \left(m e^{m+\frac{1}{r'}}\right)^p \left(1+\frac{n}{m}\right)^{\frac{p(m-1)}{r'}} \norm{u}^{mp}_{q} \norm{Q}^p_{\mathcal{P}(^m\ell^n_q)}
\end{equation}
and hence, in view of \eqref{Pal-Vasu-P5-eqn-int-01} and \eqref{Pal-Vasu-P5-eqn-004}, it follows that
\begin{equation}\label{Pal-Vasu-P5-eqn-005}
	K^p_m(B_{\ell^n_q},X,\lambda)\geq \sqrt[mp]{\lambda}\left(\pi_1(\mbox{Id}_X)\,me^{m+\frac{1}{r'}}\right)^{-1/m} \left(1+\frac{n}{m}\right)^{-\frac{m-1}{mr'}}.
\end{equation} 
On the other hand, we have 
\begin{equation*}
	\norm{u}_r\leq n^{\frac{1}{r}-\frac{1}{q}}\norm{u}_q \,\,\, \mbox{for} \,\, r<q.
\end{equation*}
Thus, we obtain
\begin{equation*}
		\sum_{\mathbf{j}\in \mathcal{J}(m,n)} \norm{x_{\mathbf{j}}(Q)}^p|u_{\mathbf{j}}|^p \leq \pi^p_1(\mbox{Id}_X) \left(m e^{m+\frac{1}{r'}}\right)^p \left(1+\frac{n}{m}\right)^{\frac{p(m-1)}{r'}} \left(n^{\frac{1}{r}-\frac{1}{q}}\right)^{p(m-1)}\norm{u}^{mp}_{q} \norm{Q}^p_{\mathcal{P}(^m\ell^n_q)}
\end{equation*}
and hence \eqref{Pal-Vasu-P5-eqn-int-01} gives
\begin{equation}\label{Pal-Vasu-P5-eqn-006}
		K^p_m(B_{\ell^n_q},X,\lambda)\geq \sqrt[mp]{\lambda}\left(\pi_1(\mbox{Id}_X)\,me^{m+\frac{1}{r'}}\right)^{-1/m} \left(n^{\frac{1}{r}-\frac{1}{q}}\right)^{-\frac{m-1}{m}}\left(1+\frac{n}{m}\right)^{-\frac{m-1}{mr'}}.
\end{equation}
From \eqref{Pal-Vasu-P5-eqn-005} and \eqref{Pal-Vasu-P5-eqn-006}, we observe that there exists a constant $E(X)=(\pi_{1}(\text{Id}_{X})m e^{m+1/r'})^{-1/m}$, depending on $X$, such that  
	$$
K^p_m(B_{\ell^n_q},X,\lambda)\geq \left\{\begin{array}{ll}
E(X)\sqrt[mp]{\lambda}\left(1+\frac{n}{m}\right)^{-\frac{m-1}{mr'}}, & \mbox{for \, $q\leq r$}\\[5mm] 
	
	E(X)\sqrt[mp]{\lambda}\left(n^{\frac{1}{r}-\frac{1}{q}}\right)^{-\frac{m-1}{m}}\left(1+\frac{n}{m}\right)^{-\frac{m-1}{mr'}}, & \mbox{for \, $r<q$}
\end{array}\right.
$$
for all $1\leq p<\infty$ and $m\in \mathbb{N}$.
The concluding part of this proof lies on the minimizing the right side of this inequality in $m$. We shall use the following important estimate
\begin{equation*}
	\left(1+\frac{n}{m}\right)^{-\frac{m-1}{m}}\geq 2^{-\frac{m-1}{m}}\min \left\{1, \left(\frac{n}{m}\right)^{-\frac{m-1}{m}}\right\}\geq \frac{1}{2}\min \left\{1, \frac{m n^{1/m}}{n m^{1/m}}\right\}.
\end{equation*}
Now, we observe that the function $h: (0, \infty)\rightarrow \mathbb{R}$ defined by $h(x)=xn^{1/x}$ attains a minimum value $e\log n$ at $x=\log n$. This implies that
\begin{equation*}
	\left(1+\frac{n}{m}\right)^{-\frac{m-1}{m}}\geq \frac{1}{2}\left(\frac{\log n}{n}\right).
\end{equation*}
Therefore, there exists a constant $C(X)$ such that
	$$
\inf_{m\in \mathbb{N}} K^p_m(B_{\ell^n_q},X,\lambda)\geq \left\{\begin{array}{ll}
	C(X)\,\left(\frac{\log n}{n}\right)^{\frac{1}{r'}}, & \mbox{for \, $q\leq r$}\\[5mm] 
	
	C(X)\,n^{\frac{1}{q}-\frac{1}{r}}\left(\frac{\log n}{n}\right)^{\frac{1}{r'}}, & \mbox{for \, $r<q$},
\end{array}\right.
$$
where $C(X)=e^{-1}\inf_{m\in \mathbb{N}}(\pi_{1}(\mbox{Id}_{X}))^{-1/m}$.
This completes the proof.
\end{pf}
 
\begin{pf}[{\bf Proof of Theorem \ref{Pal-Vasu-P5-thm-02}}]
	At first, we establish the claim of the upper bound for $ K^p(B_{\ell^n_q},X,\lambda)$. To prove this, we need an important result by Maurey and Pisier \cite[Theorem 14.5]{diestel-Abs-sum-book}: For every infinite dimensional Banach space $X$,
	\begin{equation*}
		\inf \left\{2\leq r\leq \infty : X \,\, \mbox{has cotype}\,\, r \right\}=\sup \left\{2\leq r\leq \infty : X \,\,\mbox{finitely factors}\,\, \ell_r \hookrightarrow \ell_{\infty}\right\}.
	\end{equation*}
It is known (see \cite{diestel-Abs-sum-book}) that, even when the infimum defining  $\mathrm{Cot}(X)$ is not attained, the space $X$ still finitely factors $\ell_{\mathrm{Cot}(X)} \hookrightarrow \ell_{\infty}$.	We say that $X$ finitely factors $\ell_{\text{Cot}(X)} \hookrightarrow \ell_{\infty}$ if for every $\epsilon>0$ and $n\in \mathbb{N}$, there exist $x_1, \ldots, x_n\in X$ such that 
	\begin{equation}\label{Pal-Vasu-P5-eqn-007}
		\frac{1}{1+\epsilon}\norm{z}_{\infty}\leq \norm{\sum_{j=1}^{n}z_jx_j}\leq \norm{z}_{\text{Cot}(X)}
	\end{equation}
for all $z=(z_1, \ldots, z_n)\in \mathbb{C}^n$. In particular, by considering $z=e_j=(0,\ldots, 0,1,0,\ldots,0)$, \eqref{Pal-Vasu-P5-eqn-007} reduces to
\begin{equation*}
	\frac{1}{1+\epsilon}=\frac{1}{1+\epsilon}\norm{e_j}_{\infty}\leq \norm{x_j},\,\,\, 1\leq j\leq n.
\end{equation*}
Hence, for a given $\epsilon>0$, there exist $x_1, \ldots, x_n\in X$ such that 
\begin{equation}\label{Pal-Vasu-P5-eqn-008}
	\frac{n}{(1+\epsilon)^p}\leq \sum_{j=1}^{n}\norm{x_j}^p,\,\, 1\leq p<\infty.
\end{equation}
In view of \eqref{Pal-Vasu-P5-eqn-007} and the definition of $K^p_{1}(B_{\ell^n_q},X)$ for $1$-homogeneous polynomial in $X$, for any $1\leq p<\infty$ we obtain
\begin{align}\label{Pal-Vasu-P5-eqn-009}
	\sum_{j=1}^{n}\norm{x_j}^p\frac{1}{n^{p/q}} \nonumber &\leq \frac{1}{\left(K^p_{1}(B_{\ell^n_q},X)\right)^p}\sup_{z\in B_{\ell^n_{q}}}\norm{\sum_{j=1}^{n}z_jx_j}^p \\[2mm] &\leq \frac{1}{\left(K^p_{1}(B_{\ell^n_q},X)\right)^p}\sup_{z\in B_{\ell^n_{q}}}\norm{z}^p_{\text{Cot}(X)}.
\end{align}
Now, we note that if $q\leq \text{Cot}(X)$ then $\sup_{z\in B_{\ell^n_{q}}}\norm{z}_{\text{Cot}(X)}\leq 1$. Then, using \eqref{Pal-Vasu-P5-eqn-008} and \eqref{Pal-Vasu-P5-eqn-009}, we have that
\begin{equation*}
	\left(K^p_{1}(B_{\ell^n_q}, X)\right)^p \leq \frac{(1+\epsilon)^p}{n^{1-\frac{p}{q}}}
\end{equation*}
for every $1\leq p<\infty$ and for all $\epsilon>0$. Thus, we have
\begin{equation*}
	K^p_{1}(B_{\ell^n_q}, X)\leq \frac{1}{n^{\frac{1}{p}-\frac{1}{q}}}, \,\, 1\leq p<\infty.
\end{equation*}
On the other hand, if $q>\text{Cot}(X)$ then we have 
\begin{equation*}
	\norm{z}_{\text{Cot}(X)}\leq n^{\frac{1}{\text{Cot}(X)}-\frac{1}{q}}\norm{z}_q. \quad (\mbox{ We assume the convention}\,\, n^{\frac{1}{\infty}}=1)
\end{equation*} 
In view of \eqref{Pal-Vasu-P5-eqn-008} and \eqref{Pal-Vasu-P5-eqn-009}, for any $1\leq p<\infty$ it is easy to see that
\begin{equation*}
	\left(K^p_{1}(B_{\ell^n_q}, X)\right)^p \leq \frac{(1+\epsilon)^p}{n^{1-\frac{p}{\text{Cot}(X)}}}
\end{equation*}
holds for all $\epsilon>0$. Consequently, it follows that
\begin{equation*}
	K^p_{1}(B_{\ell^n_q}, X) \leq \frac{1}{n^{\frac{1}{p}-\frac{1}{\text{Cot}(X)}}}.
\end{equation*}
Now, we have the following observation:
\begin{equation}\label{Pal-Vasu-P5-Revise-001}
	K^p(B_{\ell^n_q},X,\lambda)\leq K^p_{1}(B_{\ell^n_q},X,\lambda)=\lambda^{1/p}K^p_{1}(B_{\ell^n_q},X).
\end{equation}
Therefore, by applying \eqref{Pal-Vasu-P5-Revise-001}, we obtain 
\begin{equation}
	K^p(B_{\ell^n_q},X,\lambda)\leq 
	\begin{cases}
		\dfrac{\lambda^{1/p}}{n^{\frac{1}{p}-\frac{1}{q}}}, & \text{for } q \leq \text{Cot}(X), \\[4mm]
		\dfrac{\lambda^{1/p}}{n^{\frac{1}{p}-\frac{1}{\text{Cot}(X)}}}, & \text{for } q > \text{Cot}(X).
	\end{cases}
	\label{Pal-Vasu-P5-Revise-002}
\end{equation}

If the Banach space $X$ has no finite cotype, we have $\text{Cot}(X)=\infty$. In view of \eqref{Pal-Vasu-P5-Revise-002}, we obtain the desired upper bound for $K^p(B_{\ell^n_q},X,\lambda)$.\\ 

\vspace{2mm}
\allowdisplaybreaks
 Now we proceed for the proof of lower bound for $K^p(B_{\ell^n_q},X,\lambda)$. Suppose that $X$ has cotype $r$. To prove the lower bound of $K^p(B_{\ell^n_{q}},X, \lambda)$ in this case, we use the following technical result by Bombal \cite[Theorem 3.2]{bombal-QJM-2004}: Given an $m$-homogeneous polynomial $Q:\mathbb{C}^n\rightarrow X$, $Q(z)=\sum_{|\alpha|=m}x_{\alpha}z^{\alpha}$, we denote by $H: \mathbb{C}^n\times \cdots \times \mathbb{C}^n \rightarrow X$ the unique symmetric $m$-linear mapping associated to $Q$. Then we have
\begin{equation*}
	\left(\sum_{i\in \mathcal{M}(m,n)}\norm{H(e_{i_1}, \dots, e_{i_m})}^r\right)^{\frac{1}{r}}\leq C_r(X)^m \norm{H}_{B_{\ell^n_q}}.
\end{equation*}
Let us consider $s>1$ be such that $1/r+1/s=1$. Then, by applying H\"{o}lder's inequality, for all $z\in \mathbb{C}^n$, we obtain
\begin{align*}
	\sum_{\alpha \in\Lambda(m,n)}\norm{x_{\alpha}z^{\alpha}}^p & = \sum_{\mathbf{j}\in \mathcal{J}(m,n)} \text{card}[\mathbf{j}]^p\norm{H(e_{j_1}, \ldots, e_{j_m})}^p |z_{j_1}\cdots z_{j_m}|^p\\[2mm]& =
	\sum_{i\in \mathcal{M}(m,n)} \norm{H(e_{i_1}, \ldots, e_{i_m})}^p |z_{i_1}\cdots z_{i_m}|^p \\[2mm]&
	\leq \left(	\sum_{i\in \mathcal{M}(m,n)} \norm{H(e_{i_1}, \ldots, e_{i_m})} |z_{i_1}\cdots z_{i_m}|\right)^p\\[2mm]&
	\leq \left[\left(\sum_{i\in \mathcal{M}(m,n)}\norm{H(e_{i_1}, \ldots, e_{i_m})}^r\right)^{\frac{1}{r}}\left(\sum_{i\in \mathcal{M}(m,n)}|z_{i_1}\cdots z_{i_m}|^s\right)^{\frac{1}{s}}\right]^p\\[2mm]&
	\leq C_r(X)^{mp} \norm{H}^p_{B_{\ell^n_q}} \left(\sum_{i\in \mathcal{M}(m,n)}|z_{i_1}\cdots z_{i_m}|^s\right)^{\frac{p}{s}}\\[2mm]&
	=C_r(X)^{mp} \norm{H}^p_{B_{\ell^n_q}} \left(|z_{1}|^s+\cdots+|z_n|^s\right)^{\frac{mp}{s}}.
\end{align*}
Now, for $q\leq s$, we have $\norm{z}_s\leq \norm{z}_q$, whereas 
\begin{equation*}
	\norm{z}_s\leq n^{\frac{1}{s}-\frac{1}{q}}\norm{z}_q \,\,\, \mbox{for}\,\, s<q.
\end{equation*}
Hence, we obtain
	$$
\sum_{\alpha \in\Lambda(m,n)}\norm{x_{\alpha}z^{\alpha}}^p \leq \left\{\begin{array}{ll}
C_r(X)^{mp}\norm{H}^p_{B_{\ell^n_q}}\norm{z}^{mp}_{q}, & \mbox{for \, $q\leq s$}\\[5mm] 
	
	C_r(X)^{mp}\left(n^{\frac{1}{s}-\frac{1}{q}}\right)^{mp}\norm{H}^p_{B_{\ell^n_q}}\norm{z}^{mp}_{q}, & \mbox{for \, $q>s$}.
\end{array}\right.
$$
Therefore, by virtue of polarization formula for multilinear forms (see e.g. \cite[Section 2.6]{defant-book}), for all $z\in B_{\ell^n_q}$ we obtain 
$$
\sum_{\alpha \in\Lambda(m,n)}\norm{x_{\alpha}z^{\alpha}}^p \leq \left\{\begin{array}{ll}
	C_r(X)^{mp}\left(\frac{m^m}{m!}\right)^p\norm{Q}^p_{B_{\ell^n_q}}, & \mbox{for \, $q\leq s$}\\[5mm] 
	
	C_r(X)^{mp}\left(n^{\frac{1}{s}-\frac{1}{q}}\right)^{mp}\left(\frac{m^m}{m!}\right)^p\norm{Q}^p_{B_{\ell^n_q}}, & \mbox{for \, $q>s$}
\end{array}\right.
$$
and whence,
$$
\sum_{\alpha \in\Lambda(m,n)}\norm{x_{\alpha}z^{\alpha}}^p \leq \left\{\begin{array}{ll}
	C_r(X)^{mp}e^{mp}\norm{Q}^p_{B_{\ell^n_q}}, & \mbox{for \, $q\leq s$}\\[5mm] 
	
	C_r(X)^{mp}\left(n^{\frac{1}{s}-\frac{1}{q}}\right)^{mp}e^{mp}\norm{Q}^p_{B_{\ell^n_q}}, & \mbox{for \, $q>s$}.
\end{array}\right.
$$
In view of Definition \ref{Pal-Vasu-P5-defn-02}, for every $m\in \mathbb{N}$ and $1\leq p<\infty$ we have
$$
K^p_m(B_{\ell^n_q},X) \geq \left\{\begin{array}{ll}
	\frac{1}{e\,C_r(X)}, & \mbox{for \, $q\leq s$}\\[5mm] 
	
	\frac{1}{e\,C_r(X)n^{\frac{1}{s}-\frac{1}{q}}}, & \mbox{for \, $q>s$}.
\end{array}\right.
$$
Finally, by applying Lemma \ref{Pal-Vasu-P5-Lemma-01}(b), for every $\lambda>1$, $1\leq p<\infty$, and $1\leq q\leq \infty$, we obtain
 $$
 K^p(B_{\ell^n_q},X, \lambda) \geq \left\{\begin{array}{ll}
 	\left(\frac{\lambda-1}{\lambda}\right)^{1/p}\frac{1}{e\,C_r(X)}, & \mbox{for \, $q\leq s$}\\[5mm] 
 	
 		\left(\frac{\lambda-1}{\lambda}\right)^{1/p}\frac{1}{e\,C_r(X)n^{\frac{1}{s}-\frac{1}{q}}}, & \mbox{for \, $q>s$}.
 \end{array}\right.
 $$
 
 For Banach spaces $X$ without finite cotype the lower bound for $K^p(B_{\ell^n_q},X,\lambda)$ has been established in Corollary \ref{Pal-Vasu-P5-cor-01}. In this case, we obtain the following lower bound for $K^p(B_{\ell^n_q},X,\lambda)$:
 \begin{equation*}
 	K^p(B_{\ell^n_q},X,\lambda)\geq \left(\frac{\lambda-1}{\lambda}\right)^{\frac{1}{p}}\frac{1}{n^{\frac{1}{p}-\frac{1}{pq}}}.
 \end{equation*}
This completes the proof.
\end{pf}

\begin{pf}[{\bf Proof of Theorem \ref{Pal-Vasu-P5-thm-03}}]
	We note that the proof of (1) is just a consequence of \eqref{Pal-Vasu-P5-eqn-int-03} and Theorem \ref{Pal-Vasu-P5-thm-02}. With the help of these observations, we obtain that
	\begin{align*}
		K^p(B_{\ell^n_q},T,\lambda) & \geq \max \left\{K^p(B_{\ell^n_q},X,\lambda/\norm{T}^p), K^p(B_{\ell^n_q},Y,\lambda/\norm{T}^p)\right\} \\[2mm]&
		\geq \begin{cases}
			\left(\frac{\lambda-\norm{T}^p}{\lambda}\right)^{\frac{1}{p}}\frac{1}{e\, \min\{C_r(X),C_r(Y)\}}  & \text{if } q\leq s \\[5mm]
		 \left(\frac{\lambda-\norm{T}^p}{\lambda}\right)^{\frac{1}{p}}\frac{1}{e\, \min\{C_r(X),C_r(Y)\}}\left(\frac{1}{n}\right)^{\frac{1}{s}-\frac{1}{q}}& \text{if } q>s,
		\end{cases} 	
	\end{align*}
where $s=r/(r-1)$. Hence, the result follows.\\[1mm]

Next, for the proof of (2), we partly follow the techniques of the proof of Theorem \ref{Pal-Vasu-P5-thm-01} and hypercontractivity of vector-valued  Bohnenblust-Hille type inequality. Applying Theorem \ref{Kumar-thm-arxiv} and H\"{o}lder's inequality, there is a constant $\widetilde{C}>0$ such that for every $m$-homogeneous polynomial $Q(z)=\sum_{\alpha \in\Lambda(m,n)}x_{\alpha}z^{\alpha}\in \mathcal{P}(^m\mathbb{C}^n, X)$, we have
 \begin{align*}\label{Pal-Vasu-P5-eqn-010}
 	\sum_{\alpha \in\Lambda(m,n)}\norm{T(x_{\alpha})z^{\alpha}}^p  &\leq  \left[\left(\sum_{\alpha \in\Lambda(m,n)}1\right)^{\frac{(r-1)ms-r+s}{rms}}\times \left(\sum_{\alpha \in\Lambda(m,n)}\norm{T(x_{\alpha})z^{\alpha}}^{\frac{rsm}{r+(m-1)s}}\right)^{\frac{r+(m-1)s}{rsm}}\right]^p\\[2mm] &
 	\leq \left(\sum_{\alpha \in\Lambda(m,n)}1\right)^{\frac{((r-1)ms-r+s)p}{rms}} \widetilde{C}^{mp} \norm{Q}^p_{\ell^n_q}
 \end{align*}
the constant $\widetilde{C}$ is given by 
\begin{equation*}
	\widetilde{C}^m=\pi_{s,1}(T)M_{r}(Y)(\sqrt{2})^{m-1}m^{\frac{r-1}{r}}\left(\frac{m}{m-1}\right)^{m-1}.
\end{equation*}
We observe that
\begin{equation}\label{Pal-Vasu-P5-eqn-011}
\sum_{\alpha \in\Lambda(m,n)}1=\binom{n+m-1}{m}\leq e^m\left(1+\frac{n}{m}\right)^m.
\end{equation}
By virtue of \eqref{Pal-Vasu-P5-eqn-int-01} and \eqref{Pal-Vasu-P5-eqn-int-02} and using the fact \eqref{Pal-Vasu-P5-eqn-011}, there is a constant $E>0$ such that for every $m$, $1\leq p<\infty$, and $\norm{T}<\lambda^{1/p}$ we have
\begin{equation}\label{Pal-Vasu-P5-eqn-012}
	K^p_m(B_{\ell^n_q},T,\lambda)\geq E\, \lambda^{\frac{1}{mp}}\left(1+\frac{n}{m}\right)^{-\frac{(r-1)ms-r+s}{rms}},
\end{equation}
where $E=e\,\widetilde{C}$. By minimizing the right side of \eqref{Pal-Vasu-P5-eqn-012}, as demonstrated in the proof of Theorem \ref{Pal-Vasu-P5-thm-01}, we observe that there is a constant $D>0$ such that for $1\leq p<\infty$ and $\norm{T}<\lambda^{1/p}$, we obtain
\begin{equation*}
	\inf_{m\in \mathbb{N}} \left\{K^p_m(B_{\ell^n_q},T,\lambda)\right\}\geq D\, \left(\frac{\log n}{n}\right)^{1-\frac{1}{r}},
\end{equation*}
where $D=(1/2)^{1-1/r}e\widetilde{C}$.
Finally, using Lemma \ref{Pal-Vasu-P5-Lemma-01}(a), we obtain desired conclusion. This completes the proof.
\end{pf}

\noindent\textbf{Acknowledgment:} 
The authors would like to express their sincerest gratitude to the referee for careful reading of the manuscript and many valuable suggestions, which greatly helped to improve the clarity of the exposition in this manuscript. 
The research of the first named author is supported by NBHM Research Grant (No. 02011/19/2025/NBHM (R.P.)/R\&D II/9461), DAE, Govt. of India.
The research of the second named author is supported by DST-INSPIRE Fellowship (IF 190721), Department of Science and Technology, New Delhi, India.

\end{document}